\documentclass[numbers=enddot,12pt,final,onecolumn,notitlepage]{scrartcl}%
\usepackage[headsepline,footsepline,manualmark]{scrlayer-scrpage}
\usepackage[all,cmtip]{xy}
\usepackage{amssymb}
\usepackage{amsmath}
\usepackage{amsthm}
\usepackage{framed}
\usepackage{comment}
\usepackage{color}
\usepackage[sc]{mathpazo}
\usepackage[T1]{fontenc}
\usepackage{needspace}
\usepackage{tabls}
\usepackage[breaklinks=true]{hyperref}
\providecommand{\U}[1]{\protect\rule{.1in}{.1in}}
\newcounter{exer}

\numberwithin{exer}{section}
\theoremstyle{definition}
\newtheorem{theo}{Theorem}[section]
\newenvironment{theorem}[1][]
{\begin{theo}[#1]\begin{leftbar}}
{\end{leftbar}\end{theo}}
\newtheorem{lem}[theo]{Lemma}
\newenvironment{lemma}[1][]
{\begin{lem}[#1]\begin{leftbar}}
{\end{leftbar}\end{lem}}
\newtheorem{prop}[theo]{Proposition}

\newtheorem{defi}[theo]{Definition}
\newenvironment{definition}[1][]
{\begin{defi}[#1]\begin{leftbar}}
{\end{leftbar}\end{defi}}
\newtheorem{remk}[theo]{Remark}
\newenvironment{remark}[1][]
{\begin{remk}[#1]\begin{leftbar}}
{\end{leftbar}\end{remk}}
\newtheorem{coro}[theo]{Corollary}

\newtheorem{conv}[theo]{Convention}

\newtheorem{quest}[theo]{Question}

\newtheorem{warn}[theo]{Warning}

\newtheorem{conj}[theo]{Conjecture}

\newtheorem{exam}[theo]{Example}

\newtheorem{exmp}[exer]{Exercise}

\newenvironment{statement}{\begin{quote}}{\end{quote}}
\let\sumnonlimits\sum
\let\prodnonlimits\prod
\let\cupnonlimits\bigcup
\let\capnonlimits\bigcap
\renewcommand{\sum}{\sumnonlimits\limits}
\renewcommand{\prod}{\prodnonlimits\limits}
\renewcommand{\bigcup}{\cupnonlimits\limits}
\renewcommand{\bigcap}{\capnonlimits\limits}
\excludecomment{verlong}
\includecomment{vershort}
\excludecomment{noncompile}

\ihead{Zeckendorf family identities generalized}
\ohead{page \thepage}
\begin{document}

\title{Zeckendorf family identities generalized}
\author{Darij Grinberg}
\date{
April 13, 2026, brief version}
\maketitle

\begin{abstract}
\textbf{Abstract.} In \cite{1}, Philip Matchett Wood and Doron Zeilberger have
constructed identities for the Fibonacci numbers $f_{n}$ of the form%
\begin{align*}
1f_{n}  &  =f_{n}\text{ for all }n\geq1;\\
2f_{n}  &  =f_{n-2}+f_{n+1}\text{ for all }n\geq3;\\
3f_{n}  &  =f_{n-2}+f_{n+2}\text{ for all }n\geq3;\\
4f_{n}  &  =f_{n-2}+f_{n}+f_{n+2}\text{ for all }n\geq3;\\
&  \text{etc.;}\\
kf_{n}  &  =\sum_{s\in S_{k}}f_{n+s}\text{ for all }n>\max\left\{  -s\mid s\in
S_{k}\right\}  \text{,}%
\end{align*}
where $S_{k}$ is a fixed ``lacunar'' set of integers (``lacunar'' means that
no two elements of this set are consecutive integers) depending only on $k$.
(The condition $n>\max\left\{  -s\mid s\in S_{k}\right\}  $ is only to make
sure that all addends $f_{n+s}$ are well-defined. If the Fibonacci sequence is
properly continued to the negative, this condition drops out.)\newline In this
note we prove a generalization of these identities: For any family $\left(
a_{1},a_{2},\ldots,a_{p}\right)  $ of integers, there exists one and only one
finite lacunar set $S$ of integers such that every $n$ (high enough to make
the Fibonacci numbers in the equation below well-defined) satisfies
\[
f_{n+a_{1}}+f_{n+a_{2}}+\cdots+f_{n+a_{p}}=\sum\limits_{s\in S}f_{n+s}.
\]
The proof uses the Fibonacci-approximating properties of the golden ratio. It
would be interesting to find a purely combinatorial proof.

\end{abstract}


\hrule

\begin{statement}
This is a brief version of my note \cite{this.long}. For a long version, which
gives more details in the proofs, see \cite{this.long}.
\end{statement}

\hrule

\section{The main result}

The purpose of this note is to establish a generalization of the so-called
\textit{Zeckendorf family identities} which were discussed in \cite{1} and
\cite{Gerdem08}. Before we can state it, we need a few definitions:

\begin{definition}
A subset $S$ of $\mathbb{Z}$ is called \textit{lacunar} if it satisfies
$\left(  s+1\notin S\text{ for every }s\in S\right)  $.
\end{definition}

In other words, a subset $S$ of $\mathbb{Z}$ is lacunar if and only if it
contains no two consecutive integers.

\begin{definition}
\label{def.fib}The \textit{Fibonacci sequence} $\left(  f_{1},f_{2}%
,f_{3},\ldots\right)  $ is a sequence of positive integers defined recursively
by the initial values $f_{1}=1$ and $f_{2}=1$ and the recurrence relation
$\left(  f_{n}=f_{n-1}+f_{n-2}\text{ for all }n\in\mathbb{N}\text{ satisfying
}n\geq3\right)  $.
\end{definition}

(Here and in the following, $\mathbb{N}$ denotes the set $\left\{
0,1,2,\ldots\right\}  $.)

\begin{remark}
Many authors define the Fibonacci sequence slightly differently: They define
it as a sequence $\left(  f_{0},f_{1},f_{2},\ldots\right)  $ of nonnegative
integers defined recursively by the initial values $f_{0}=0$ and $f_{1}=1$ and
the recurrence relation $\left(  f_{n}=f_{n-1}+f_{n-2}\text{ for all }%
n\in\mathbb{N}\text{ satisfying }n\geq2\right)  $. Thus, this sequence begins
with a $0$, unlike the Fibonacci sequence defined in our Definition
\ref{def.fib}. However, starting at its second term $f_{1}=1$, this sequence
takes precisely the same values as the Fibonacci sequence defined in our
Definition \ref{def.fib} (because both sequences satisfy $f_{1}=1$ and
$f_{2}=1$, and from here on the recurrence relation ensures that their values
remain equal). So it differs from the latter sequence only in the presence of
one extra term $f_{0}=0$ at the front.
\end{remark}

The Fibonacci sequence is one of the best known integer sequences from
combinatorics. It has had conferences, books and a journal devoted to it. By
way of example, let us only mention Vorobiev's book \cite{Vorobi02}, which is
entirely concerned with Fibonacci numbers, and Benjamin's and Quinn's text
\cite{BenQui03} on bijective proofs, which includes many identities for
Fibonacci numbers.

In \cite{1}, Wood and Zeilberger discuss bijective proofs of the so-called
\textit{Zeckendorf family identities} (see also  \cite{Gerdem08}). These
identities are a family of identities for Fibonacci numbers (one for each
positive integer); the first seven of these identities are%
\begin{align*}
1f_{n} &  =f_{n}\text{ for all }n\geq1;\\
2f_{n} &  =f_{n-2}+f_{n+1}\text{ for all }n\geq3;\\
3f_{n} &  =f_{n-2}+f_{n+2}\text{ for all }n\geq3;\\
4f_{n} &  =f_{n-2}+f_{n}+f_{n+2}\text{ for all }n\geq3;\\
5f_{n} &  =f_{n-4}+f_{n-1}+f_{n+3}\text{ for all }n\geq5;\\
6f_{n} &  =f_{n-4}+f_{n+1}+f_{n+3}\text{ for all }n\geq5;\\
7f_{n} &  =f_{n-4}+f_{n+4}\text{ for all }n\geq5.
\end{align*}
In general, for each positive integer $k$, the $k$-th Zeckendorf family
identity expresses $kf_{n}$ (for each sufficiently large integer $n$) as a sum
of the form $\sum\limits_{s\in S}f_{n+s}$, where $S$ is a finite lacunar
subset of $\mathbb{Z}$. Of course, the subset $S$ does not depend on $n$.

Our main theorem is the following:

\begin{theorem}
[generalized Zeckendorf family identities]\label{thm.1} Let $T$ be a finite
set, and let $a_{t}$ be an integer for every $t\in T$.

Then, there exists one and only one finite lacunar subset $S$ of $\mathbb{Z}$
such that\footnotemark%
\[
\left(
\begin{array}
[c]{c}%
\sum\limits_{t\in T}f_{n+a_{t}}=\sum\limits_{s\in S}f_{n+s}\text{ for every
}n\in\mathbb{Z}\text{ which}\\
\text{satisfies }n>\max\left(  \left\{  -a_{t}\mid t\in T\right\}
\cup\left\{  -s\mid s\in S\right\}  \right)
\end{array}
\right)  .
\]

\end{theorem}

\footnotetext{Here and in the following, $\max\varnothing$ is understood to be
$0$.}

\begin{remark}
\begin{enumerate}
\item The \textit{Zeckendorf family identities}\ from \cite{1} are the result
of applying Theorem~\ref{thm.1} to the case when all $a_{t}$ are $=0$.

\item The condition $n>\max\left(  \left\{  -a_{t}\mid t\in T\right\}
\cup\left\{  -s\mid s\in S\right\}  \right)  $ in Theorem~\ref{thm.1} serves
only to ensure that the Fibonacci numbers $f_{n+a_{t}}$ for all $t\in T$ and
$f_{n+s}$ for all $s\in S$ are well-defined. (If we would define the Fibonacci
numbers $f_{n}$ for integers $n\leq0$ by extending the recurrence relation
$f_{n}=f_{n-1}+f_{n-2}$ \textquotedblleft to the left\textquotedblright, then
we could drop this condition.)
\end{enumerate}
\end{remark}

The proof we shall give for Theorem~\ref{thm.1} is not combinatorial. It will
use inequalities and the properties of the golden ratio; in a sense, its
underlying ideas come from analysis (although it will not actually use any
results from analysis).

\section{Basics on the Fibonacci sequence}

We begin with some lemmas and notations:

We denote by $\mathbb{N}$ the set $\left\{  0,1,2,\ldots\right\}  $. Also, we
denote by $\mathbb{N}_{2}$ the set $\left\{  2,3,4,\ldots\right\}
=\mathbb{N}\setminus\left\{  0,1\right\}  $.

Also, let
\[
\phi=\dfrac{1+\sqrt{5}}{2}\ \ \ \ \ \ \ \ \ \ \text{and}%
\ \ \ \ \ \ \ \ \ \ \psi=\dfrac{1-\sqrt{5}}{2}.
\]
These numbers $\phi$ and $\psi$ are known as the \textit{golden ratios}. We
notice that $\phi\approx1.618\ldots$ and $\phi^{2}=\phi+1$ and $\psi
\approx-0.618\ldots$ and $\psi^{2}=\psi+1$ and $\phi+\psi=1$ and $\phi
-\psi=\sqrt{5}$.

We recall some basic and well-known facts about the Fibonacci sequence:

\begin{lemma}
\label{lem.2} Let $S$ be a finite lacunar subset of $\mathbb{N}_{2}$. Then,
$\sum\limits_{t\in S}f_{t}<f_{\max S+1}$.
\end{lemma}

\begin{proof}
WLOG assume that the set $S$ is nonempty (else, the lemma follows from our
convention that $\max\varnothing=0$). Write the set $S$ in the form $\left\{
s_{1},s_{2},\ldots,s_{k}\right\}  $ with $s_{1}<s_{2}<\cdots<s_{k}$. Every
$i\in\left\{  1,2,\ldots,k-1\right\}  $ satisfies $s_{i}+1\leq s_{i+1}-1$
(because the set $S$ is lacunar, so $s_{i+1}$ cannot be $s_{i}+1$, whence
$s_{i+1}>s_{i}+1$ and thus $s_{i+1}-1\geq s_{i}+1$), so that
\begin{equation}
f_{s_{i}+1}\leq f_{s_{i+1}-1} \label{pf.lem.2.leq}%
\end{equation}
(since the Fibonacci sequence $\left(  f_{1},f_{2},f_{3},\ldots\right)  $ is
weakly increasing). Thus,%
\begin{align*}
\sum\limits_{t\in S}f_{t}  &  =\sum\limits_{i=1}^{k}\underbrace{f_{s_{i}}%
}_{\substack{=f_{s_{i}+1}-f_{s_{i}-1}\\\text{(since }f_{s_{i}+1}=f_{s_{i}%
}+f_{s_{i}-1}\text{)}}}=\sum\limits_{i=1}^{k}\left(  f_{s_{i}+1}-f_{s_{i}%
-1}\right)  =\underbrace{\sum\limits_{i=1}^{k}f_{s_{i}+1}}_{=\sum
\limits_{i=1}^{k-1}f_{s_{i}+1}+f_{s_{k}+1}}-\underbrace{\sum\limits_{i=1}%
^{k}f_{s_{i}-1}}_{=f_{s_{1}-1}+\sum\limits_{i=2}^{k}f_{s_{i}-1}}\\
&  =\left(  \sum\limits_{i=1}^{k-1}\underbrace{f_{s_{i}+1}}_{\substack{\leq
f_{s_{i+1}-1}\\\text{(by (\ref{pf.lem.2.leq}))}}}+\,f_{s_{k}+1}\right)
-\left(  f_{s_{1}-1}+\sum\limits_{i=2}^{k}f_{s_{i}-1}\right) \\
&  \leq\left(  \sum\limits_{i=1}^{k-1}f_{s_{i+1}-1}+f_{s_{k}+1}\right)
-\left(  f_{s_{1}-1}+\sum\limits_{i=2}^{k}f_{s_{i}-1}\right) \\
&  =\left(  \sum\limits_{i=2}^{k}f_{s_{i}-1}+f_{s_{k}+1}\right)  -\left(
f_{s_{1}-1}+\sum\limits_{i=2}^{k}f_{s_{i}-1}\right) \\
&  \ \ \ \ \ \ \ \ \ \ \left(  \text{here, we substituted }i\text{ for
}i+1\text{ in the first sum}\right) \\
&  =f_{s_{k}+1}-f_{s_{1}-1}<f_{s_{k}+1}\ \ \ \ \ \ \ \ \ \ \left(  \text{since
}f_{s_{1}-1}>0\right) \\
&  =f_{\max S+1}%
\end{align*}
(since $s_{k}=\max S$), which proves Lemma~\ref{lem.2}. (An alternative proof
proceeds by strong induction over $\max S$; it uses $f_{\max S+1}=f_{\max
S}+f_{\max S-1}$ in the induction step.)
\end{proof}

\begin{lemma}
[existence part of the Zeckendorf theorem]\label{lem.3} Let $n \in\mathbb{N}$.
Then, there exists a finite lacunar subset $T$ of $\mathbb{N}_{2}$ such that
$n=\sum\limits_{t\in T}f_{t}$.
\end{lemma}

\begin{proof}
Strong induction over $n$. The case $n=0$ needs to be treated separately. In
the induction step for $n>0$, the main idea is to let $t_{1}$ be the maximal
$\tau\in\mathbb{N}_{2}$ satisfying $f_{\tau}\leq n$ (this exists because
$f_{2}=1\leq n$ and because the Fibonacci sequence is increasing and unbounded
from above), and to apply Lemma~\ref{lem.3} to $n-f_{t_{1}}$ instead of $n$.
(This yields a finite lacunar subset $T^{\prime}$ of $\mathbb{N}_{2}$
satisfying $n-f_{t_{1}}=\sum_{t\in T^{\prime}}f_{t}$; now, it remains to be
shown that the set $T^{\prime}\cup\left\{  t_{1}\right\}  $ is still lacunar.
To check this, observe that $n<f_{t_{1}+1}$, so that $n-f_{t_{1}}<f_{t_{1}%
+1}-f_{t_{1}}=f_{t_{1}-1}$, which shows that no addend $f_{t}$ of the sum
$\sum_{t\in T^{\prime}}f_{t}$ can be $f_{t_{1}-1}$ or larger.) The details are
left to the reader (and can be found in \cite{this.long}).
\end{proof}

\begin{lemma}
[uniqueness part of the Zeckendorf theorem]\label{lem.4} Let $n\in\mathbb{N}$,
and let $T$ and $T^{\prime}$ be two finite lacunar subsets of $\mathbb{N}_{2}$
such that $n=\sum\limits_{t\in T}f_{t}$ and $n=\sum\limits_{t\in T^{\prime}%
}f_{t}$. Then, $T=T^{\prime}$.
\end{lemma}

\begin{proof}
Strong induction over $n$. In the induction step for $n>0$, use
Lemma~\ref{lem.2} to show that $\max T<\max T^{\prime}+1$ and $\max T^{\prime
}<\max T+1$; these together result in $\max T=\max T^{\prime}$. Hence, the
sets $T$ and $T^{\prime}$ have an element in common, and we can reduce the
situation to one with a smaller $n$ by removing this common element from both sets.
\end{proof}

Lemmata~\ref{lem.3} and~\ref{lem.4} together yield the following theorem:

\begin{theorem}
[Zeckendorf theorem]\label{thm.5} Let $n \in\mathbb{N}$. Then, there exists
one and only one finite lacunar subset $T$ of $\mathbb{N}_{2}$ such that
$n=\sum\limits_{t\in T}f_{t}$.
\end{theorem}

Theorem~\ref{thm.5} is a classical result that can be found in various places
(e.g., \cite{Hender16}). Hoggatt proved a generalization of
Theorem~\ref{thm.5} in \cite{Hoggat72}.

\begin{definition}
Let $n\in\mathbb{N}$. Theorem~\ref{thm.5} shows that there exists one and only
one finite lacunar subset $T$ of $\mathbb{N}_{2}$ such that $n=\sum
\limits_{t\in T}f_{t}$. We will denote this set $T$ by $Z_{n}$. Thus,
$n=\sum\limits_{t\in Z_{n}}f_{t}$.
\end{definition}

\section{Inequalities for the golden ratio}

Next, we state a completely straightforward (and well-known, cf. \cite[Chapter
9, Corollary 34]{BenQui03}) theorem, which shows that the Fibonacci sequence
grows roughly exponentially (with the exponent being the golden ratio $\phi$):

\begin{theorem}
\label{thm.6} For every positive integer $n$, we have $f_{n+1}-\phi f_{n}%
=\psi^{n}$.
\end{theorem}

\begin{proof}
Binet's formula (see, e.g., \cite[Identity 240]{BenQui03} or \cite[(1.20)]%
{Vorobi02}) yields $f_{n}=\dfrac{\phi^{n}-\psi^{n}}{\sqrt{5}}$ and
$f_{n+1}=\dfrac{\phi^{n+1}-\psi^{n+1}}{\sqrt{5}}$; the rest is computation.
\end{proof}

Let us show yet another lemma for later use:

\begin{lemma}
\label{lem.7} Let $S$ be a finite lacunar subset of $\mathbb{N}_{2}$. Then,
$\sum\limits_{s\in S}\left\vert \psi\right\vert ^{s}\leq\left\vert
\psi\right\vert $.
\end{lemma}

\begin{proof}
[Proof of Lemma~\ref{lem.7}.]We note that $\psi=1-\phi$ (since $\phi+\psi=1$),
thus $\psi<0$ and therefore $\left\vert \psi\right\vert =-\psi=-\left(
1-\phi\right)  =\phi-1$.

Since $S$ is a lacunar subset of $\mathbb{N}_{2}$, the smallest element of $S$
is at least $2$, the second smallest element of $S$ is at least $4$ (since it
is larger than the smallest element by at least $2$), the third smallest
element of $S$ is at least $6$ (since it is larger than the second smallest
element by at least $2$), and so on. Since the function $\mathbb{N}%
\rightarrow\mathbb{R}$, $s\mapsto\left\vert \psi\right\vert ^{s}$ is weakly
decreasing (as $0\leq\left\vert \psi\right\vert \leq1$), we thus have%
\[
\sum_{s\in S}\left\vert \psi\right\vert ^{s}\leq\sum_{s\in\left\{
2,4,6,\ldots\right\}  }\left\vert \psi\right\vert ^{s}=\sum_{t\in\left\{
1,2,3,\ldots\right\}  }\left\vert \psi\right\vert ^{2t}=\left\vert
\psi\right\vert ^{2}\cdot\dfrac{1}{1-\left\vert \psi\right\vert ^{2}}%
\]
(by the formula for the sum of the geometric series). Since $\left\vert
\psi\right\vert ^{2}=\psi^{2}$, this rewrites as
\begin{align*}
\sum_{s\in S}\left\vert \psi\right\vert ^{s}  &  \leq\psi^{2}\cdot\dfrac
{1}{1-\psi^{2}}=\psi^{2}\cdot\dfrac{1}{1-\left(  1+\psi\right)  }%
\ \ \ \ \ \ \ \ \ \ \left(  \text{since }\psi^{2}=1+\psi\right) \\
&  =\psi^{2}\cdot\dfrac{1}{-\psi}=-\psi=\left\vert \psi\right\vert .
\end{align*}
This proves Lemma~\ref{lem.7}.
\end{proof}

\section{Proving Theorem~\ref{thm.1}}

Let us now come to the proof of Theorem~\ref{thm.1}. First, we formulate the
existence part of this theorem:

\begin{theorem}
[existence part of the generalized Zeckendorf family identities]\label{thm.8}
Let $T$ be a finite set, and let $a_{t}$ be an integer for every $t\in T$.

Then, there exists a finite lacunar subset $S$ of $\mathbb{Z}$ such that
\[
\left(
\begin{array}
[c]{c}%
\sum\limits_{t\in T}f_{n+a_{t}}=\sum\limits_{s\in S}f_{n+s}\text{ for every
}n\in\mathbb{Z}\text{ which}\\
\text{satisfies }n>\max\left(  \left\{  -a_{t}\mid t\in T\right\}
\cup\left\{  -s\mid s\in S\right\}  \right)
\end{array}
\right)  .
\]

\end{theorem}

Before we start proving this, let us introduce a notation:

\begin{definition}
Let $K$ be a subset of $\mathbb{Z}$, and let $a\in\mathbb{Z}$. Then, $K+a$
will denote the subset $\left\{  k+a\ \mid\ k\in K\right\}  $ of $\mathbb{Z}$.
\end{definition}

Clearly, $\left(  K+a\right)  +b=K+\left(  a+b\right)  $ for any two integers
$a$ and $b$. Also, $K+0=K$. Finally, if $K$ is a lacunar subset of
$\mathbb{Z}$, and if $a\in\mathbb{Z}$, then $K+a$ is lacunar as well.

\begin{proof}
[Proof of Theorem~\ref{thm.8}.]Choose a high enough integer $N$. Here,
\textquotedblleft high enough\textquotedblright\ means that $N$ should satisfy
$N\in\mathbb{N}_{2}$ and $N>\max\left\{  -a_{t}\mid t\in T\right\}  $ and%
\begin{equation}
\left\vert \psi\right\vert ^{N}\sum\limits_{t\in T}\left\vert \psi\right\vert
^{a_{t}}+\left\vert \psi\right\vert <1. \label{pf.thm.8.how-high-should-N-be}%
\end{equation}
(Such an $N$ can indeed be found\footnote{\textit{Proof.} We have $\left\vert
\psi\right\vert ^{N}\rightarrow0$ for $N\rightarrow\infty$ (since
$0<\left\vert \psi\right\vert <1$). Therefore, the left hand side of
(\ref{pf.thm.8.how-high-should-N-be}) tends to $\left\vert \psi\right\vert $
as $N\rightarrow\infty$. Thus, for all sufficiently high $N$, the left hand
side of (\ref{pf.thm.8.how-high-should-N-be}) will be $<1$, because
$\left\vert \psi\right\vert <1$. So, if we take $N$ sufficiently high, then
(\ref{pf.thm.8.how-high-should-N-be}) will hold. Of course, our other two
requirements on $N$ (namely, $N\in\mathbb{N}_{2}$ and $N>\max\left\{
-a_{t}\mid t\in T\right\}  $) can also be achieved by taking $N$ sufficiently
high.}.)

Let $\nu=\sum\limits_{t\in T}f_{N+a_{t}}$. Then, $Z_{\nu}$ is a finite lacunar
subset of $\mathbb{N}_{2}$ satisfying $\nu=\sum\limits_{t\in Z_{\nu}}f_{t}$.
Hence, Lemma~\ref{lem.7} yields
\begin{equation}
\sum\limits_{s\in Z_{\nu}}\left\vert \psi\right\vert ^{s}\leq\left\vert
\psi\right\vert . \label{pf.thm.8.lacZnu}%
\end{equation}

Define a subset $S$ of $\mathbb{Z}$ by $S=Z_{\nu}+\left(  -N\right)  $. Then,
$S$ is a finite lacunar subset of $\mathbb{Z}$ (since $Z_{\nu}$ is a finite
lacunar subset of $\mathbb{Z}$). Furthermore, from $S=Z_{\nu}+\left(
-N\right)  $, we obtain $Z_{\nu}=S+N$. Thus, the map $S\rightarrow Z_{\nu
},\ s\mapsto N+s$ is a bijection. This allows us to substitute $N+s$ for $t$
in sums over all $t\in Z_{\nu}$; we thus obtain%
\begin{align}
\sum\limits_{t\in Z_{\nu}}f_{t}  &  =\sum\limits_{s\in S}f_{N+s}%
\ \ \ \ \ \ \ \ \ \ \text{and}\nonumber\\
\sum\limits_{t\in Z_{\nu}}\left\vert \psi\right\vert ^{t}  &  =\sum
\limits_{s\in S}\left\vert \psi\right\vert ^{N+s}. \label{pf.thm.8.subs2}%
\end{align}
Hence,
\begin{equation}
\sum\limits_{t\in T}f_{N+a_{t}}=\nu=\sum\limits_{t\in Z_{\nu}}f_{t}%
=\sum\limits_{s\in S}f_{N+s}, \label{pf.thm.8.subs1a}%
\end{equation}
while the equality (\ref{pf.thm.8.subs2}) yields%
\begin{equation}
\sum\limits_{s\in S}\left\vert \psi\right\vert ^{N+s}=\sum\limits_{t\in
Z_{\nu}}\left\vert \psi\right\vert ^{t}=\sum\limits_{s\in Z_{\nu}}\left\vert
\psi\right\vert ^{s}\leq\left\vert \psi\right\vert \label{pf.thm.8.lacZnu2}%
\end{equation}
(by (\ref{pf.thm.8.lacZnu})).

So, we have chosen a high $N$ and found a finite lacunar subset $S$ of
$\mathbb{Z}$ which satisfies $\sum\limits_{t\in T}f_{N+a_{t}}=\sum
\limits_{s\in S}f_{N+s}$. But Theorem~\ref{thm.8} is not proven yet:
Theorem~\ref{thm.8} requires us to prove that there exists \textit{one} finite
lacunar subset $S$ of $\mathbb{Z}$ which works for \textit{every} $N$, while
at the moment we cannot be sure yet whether different $N$'s wouldn't produce
\textit{different} sets $S$. And, in fact, different $N$'s \textit{can}
produce different sets $S$, but (fortunately!) only if the $N$'s are too
small. As we have taken $N$ high enough, the set $S$ that we obtained turns
out to be \textit{universal}, i.e., it satisfies
\begin{equation}
\left(
\begin{array}
[c]{c}%
\sum\limits_{t\in T}f_{n+a_{t}}=\sum\limits_{s\in S}f_{n+s}\text{ for every
}n\in\mathbb{Z}\text{ which}\\
\text{satisfies }n>\max\left(  \left\{  -a_{t}\mid t\in T\right\}
\cup\left\{  -s\mid s\in S\right\}  \right)
\end{array}
\right)  . \label{BigLemma}%
\end{equation}
We are now going to prove this.

In order to prove (\ref{BigLemma}), we need two assertions:

\begin{statement}
\textit{Assertion 1:} If some $n\in\mathbb{Z}$ satisfies $n\geq N$ and
$\sum\limits_{t\in T}f_{n+a_{t}}=\sum\limits_{s\in S}f_{n+s}$, then
$\sum\limits_{t\in T}f_{\left(  n+1\right)  +a_{t}}=\sum\limits_{s\in
S}f_{\left(  n+1\right)  +s}$.
\end{statement}

\begin{statement}
\textit{Assertion 2:} If some $n\in\mathbb{Z}$ satisfies $\sum\limits_{t\in
T}f_{n+a_{t}}=\sum\limits_{s\in S}f_{n+s}$ and $\sum\limits_{t\in T}f_{\left(
n+1\right)  +a_{t}}=\sum\limits_{s\in S}f_{\left(  n+1\right)  +s}$, then
$\sum\limits_{t\in T}f_{\left(  n-1\right)  +a_{t}}=\sum\limits_{s\in
S}f_{\left(  n-1\right)  +s}$ (if \newline$n-1>\max\left(  \left\{  -a_{t}\mid
t\in T\right\}  \cup\left\{  -s\mid s\in S\right\}  \right)  $).
\end{statement}

Obviously, Assertion 1 yields (by induction) the equality $\sum\limits_{t\in
T}f_{n+a_{t}}=\sum\limits_{s\in S}f_{n+s}$ for every $n\geq N$ (the induction
base follows from (\ref{pf.thm.8.subs1a})), and Assertion 2 then proves it for
the remaining values of $n$ (by backwards induction, or, to be more precise,
by an induction step from $n+1$ and $n$ to $n-1$). Thus, once both Assertions
1 and 2 are proven, (\ref{BigLemma}) will follow and thus Theorem~\ref{thm.8}
will be proven.

Assertion 2 follows from comparing the equalities%
\[
\sum\limits_{t\in T}\underbrace{f_{\left(  n-1\right)  +a_{t}}}%
_{\substack{=f_{n+a_{t}-1}\\=f_{n+a_{t}+1}-f_{n+a_{t}}}}=\sum\limits_{t\in
T}f_{n+a_{t}+1}-\sum\limits_{t\in T}f_{n+a_{t}}=\sum\limits_{t\in T}f_{\left(
n+1\right)  +a_{t}}-\sum\limits_{t\in T}f_{n+a_{t}}%
\]
and%
\[
\sum\limits_{s\in S}\underbrace{f_{\left(  n-1\right)  +s}}%
_{\substack{=f_{n+s-1}\\=f_{n+s+1}-f_{n+s}}}=\sum\limits_{s\in S}%
f_{n+s+1}-\sum\limits_{s\in S}f_{n+s}=\sum\limits_{s\in S}f_{\left(
n+1\right)  +s}-\sum\limits_{s\in S}f_{n+s}%
\]
(whose right hand sides are equal by the assumptions of Assertion 2); thus, it
only remains to prove Assertion 1.

So let us prove Assertion 1. Here we are going to use that $N$ is high enough
(because otherwise, Assertion 1 wouldn't hold). We have $\sum\limits_{t\in
T}f_{n+a_{t}}=\sum\limits_{s\in S}f_{n+s}$ by assumption, so that
$\sum\limits_{t\in T}f_{n+a_{t}}-\sum\limits_{s\in S}f_{n+s}=0$. Thus,%
\begin{align*}
\sum\limits_{t\in T}f_{\left(  n+1\right)  +a_{t}}-\sum\limits_{s\in
S}f_{\left(  n+1\right)  +s}  &  =\sum\limits_{t\in T}f_{\left(  n+1\right)
+a_{t}}-\sum\limits_{s\in S}f_{\left(  n+1\right)  +s}-\phi\left(
\sum\limits_{t\in T}f_{n+a_{t}}-\sum\limits_{s\in S}f_{n+s}\right) \\
&  =\sum\limits_{t\in T}\left(  f_{\left(  n+1\right)  +a_{t}}-\phi
f_{n+a_{t}}\right)  -\sum\limits_{s\in S}\left(  f_{\left(  n+1\right)
+s}-\phi f_{n+s}\right) \\
&  =\sum\limits_{t\in T}\left(  f_{n+a_{t}+1}-\phi f_{n+a_{t}}\right)
-\sum\limits_{s\in S}\left(  f_{n+s+1}-\phi f_{n+s}\right)  ,
\end{align*}
so that%
\begin{align*}
&  \left\vert \sum\limits_{t\in T}f_{\left(  n+1\right)  +a_{t}}%
-\sum\limits_{s\in S}f_{\left(  n+1\right)  +s}\right\vert =\left\vert
\sum\limits_{t\in T}\left(  f_{n+a_{t}+1}-\phi f_{n+a_{t}}\right)
-\sum\limits_{s\in S}\left(  f_{n+s+1}-\phi f_{n+s}\right)  \right\vert \\
&  \leq\sum\limits_{t\in T}\left\vert f_{n+a_{t}+1}-\phi f_{n+a_{t}%
}\right\vert +\sum\limits_{s\in S}\left\vert f_{n+s+1}-\phi f_{n+s}\right\vert
\ \ \ \ \ \ \ \ \ \ \left(  \text{by the triangle inequality}\right) \\
&  =\sum\limits_{t\in T}\underbrace{\left\vert \psi^{n+a_{t}}\right\vert
}_{=\left\vert \psi\right\vert ^{n+a_{t}}=\left\vert \psi\right\vert
^{n}\left\vert \psi\right\vert ^{a_{t}}}+\sum\limits_{s\in S}%
\underbrace{\left\vert \psi^{n+s}\right\vert }_{\substack{=\left\vert
\psi\right\vert ^{n+s}\leq\left\vert \psi\right\vert ^{N+s}\\\text{(since
}n\geq N\text{ and }0<\left\vert \psi\right\vert <1\text{)}}%
}\ \ \ \ \ \ \ \ \ \ \left(  \text{by Theorem~\ref{thm.6}}\right) \\
&  \leq\sum\limits_{t\in T}\underbrace{\left\vert \psi\right\vert ^{n}%
}_{\substack{\leq\left\vert \psi\right\vert ^{N}\\\text{(since }n\geq N\text{
and }0<\left\vert \psi\right\vert <1\text{)}}}\left\vert \psi\right\vert
^{a_{t}}+\underbrace{\sum\limits_{s\in S}\left\vert \psi\right\vert ^{N+s}%
}_{\substack{\leq\left\vert \psi\right\vert \\\text{(by
(\ref{pf.thm.8.lacZnu2}))}}}\\
&  \leq\sum\limits_{t\in T}\left\vert \psi\right\vert ^{N}\left\vert
\psi\right\vert ^{a_{t}}+\left\vert \psi\right\vert =\left\vert \psi
\right\vert ^{N}\sum\limits_{t\in T}\left\vert \psi\right\vert ^{a_{t}%
}+\left\vert \psi\right\vert <1
\end{align*}
(by (\ref{pf.thm.8.how-high-should-N-be})). This leads to $\left\vert
\sum\limits_{t\in T}f_{\left(  n+1\right)  +a_{t}}-\sum\limits_{s\in
S}f_{\left(  n+1\right)  +s}\right\vert =0$ (since $\left\vert \sum
\limits_{t\in T}f_{\left(  n+1\right)  +a_{t}}-\sum\limits_{s\in S}f_{\left(
n+1\right)  +s}\right\vert $ is a nonnegative integer). In other words,
$\sum\limits_{t\in T}f_{\left(  n+1\right)  +a_{t}}=\sum\limits_{s\in
S}f_{\left(  n+1\right)  +s}$. This completes the proof of Assertion 1, and,
with it, the proof of Theorem~\ref{thm.8}.
\end{proof}

All that remains now is the (rather trivial) uniqueness part of
Theorem~\ref{thm.1}:

\begin{lemma}
[uniqueness part of the generalized Zeckendorf family identities]\label{lem.9}
Let $T$ be a finite set, and let $a_{t}$ be an integer for every $t\in T$.

Let $S$ be a finite lacunar subset of $\mathbb{Z}$ such that%
\begin{equation}
\left(
\begin{array}
[c]{c}%
\sum\limits_{t\in T}f_{n+a_{t}}=\sum\limits_{s\in S}f_{n+s}\text{ for every
}n\in\mathbb{Z}\text{ which}\\
\text{satisfies }n>\max\left(  \left\{  -a_{t}\mid t\in T\right\}
\cup\left\{  -s\mid s\in S\right\}  \right)
\end{array}
\right)  . \label{eq.lem.9.ass1}%
\end{equation}
Let $S^{\prime}$ be a finite lacunar subset of $\mathbb{Z}$ such that%
\begin{equation}
\left(
\begin{array}
[c]{c}%
\sum\limits_{t\in T}f_{n+a_{t}}=\sum\limits_{s\in S^{\prime}}f_{n+s}\text{ for
every }n\in\mathbb{Z}\text{ which}\\
\text{satisfies }n>\max\left(  \left\{  -a_{t}\mid t\in T\right\}
\cup\left\{  -s\mid s\in S^{\prime}\right\}  \right)
\end{array}
\right)  . \label{eq.lem.9.ass2}%
\end{equation}
Then, $S=S^{\prime}$.
\end{lemma}

\begin{proof}
[Proof of Lemma~\ref{lem.9}.]Let%
\begin{equation}
n=\max\left(  \left\{  -a_{t}\mid t\in T\right\}  \cup\left\{  -s\mid s\in
S\right\}  \cup\left\{  -s\mid s\in S^{\prime}\right\}  \right)  +2.
\label{l7}%
\end{equation}
Then, $n$ satisfies $n>\max\left(  \left\{  -a_{t}\mid t\in T\right\}
\cup\left\{  -s\mid s\in S\right\}  \right)  $. Thus, (\ref{eq.lem.9.ass1})
yields
\[
\sum\limits_{t\in T}f_{n+a_{t}}=\sum\limits_{s\in S}f_{n+s}=\sum\limits_{t\in
S+n}f_{t}%
\]
(here, we substituted $t$ for $n+s$, since the map $S\rightarrow
S+n,\ s\mapsto n+s$ is a bijection). Similarly, $\sum\limits_{t\in
T}f_{n+a_{t}}=\sum\limits_{t\in S^{\prime}+n}f_{t}$. Since the sets $S+n$ and
$S^{\prime}+n$ are both lacunar (since so are $S$ and $S^{\prime}$) and finite
(since so are $S$ and $S^{\prime}$), and are subsets of $\mathbb{N}_{2}$ (by
(\ref{l7})), we can now conclude from Lemma~\ref{lem.4} (applied to
$\sum\limits_{t\in T}f_{n+a_{t}}$, $S+n$ and $S^{\prime}+n$ instead of $n$,
$S$ and $S^{\prime}$) that $S+n=S^{\prime}+n$, so that $S=S^{\prime}$. This
proves Lemma~\ref{lem.9}.
\end{proof}

\begin{proof}
[Proof of Theorem~\ref{thm.1}.]Now, Theorem~\ref{thm.1} is clear, since the
existence follows from Theorem~\ref{thm.8} and the uniqueness from
Lemma~\ref{lem.9}.
\end{proof}


\begin{thebibliography}{99999999}                                                                                         %


\bibitem[BenQui03]{BenQui03}Arthur T. Benjamin and Jennifer J. Quinn,
\textit{Proofs that Really Count: The Art of Combinatorial Proof}, The
Mathematical Association of America, 2003.

\bibitem[Gerdem08]{Gerdem08}%
\href{https://www.fq.math.ca/Papers1/46_47-3/Gerdemann2.pdf}{Dale Gerdemann,
\textit{Combinatorial proofs of Zeckendorf family identities}, The Fibonacci
Quarterly \textbf{46/47} (2008/09), no. 3, pp. 249--261.}

\bibitem[Grinbe11]{this.long}Darij Grinberg, \textit{Zeckendorf family
identities generalized},
April 13, 2026, *long version*.\newline%
\url{http://www.cip.ifi.lmu.de/~grinberg/zeckendorfLONG.pdf}\newline Also
available as an ancillary file to
\href{https://arxiv.org/abs/1103.4507v3}{arXiv preprint arXiv:1103.4507v3}.

\bibitem[Hender16]{Hender16}Nik Henderson, \textit{What is Zeckendorf's
Theorem?}, 23 July 2016.\newline\url{https://math.osu.edu/sites/math.osu.edu/files/henderson_zeckendorf.pdf}

\bibitem[Hoggat72]{Hoggat72}V. E. Hoggatt, Jr., \textit{Generalized Zeckendorf
theorem}, The Fibonacci Quarterly \textbf{10} (1972), Issue 1, pp.
89--94.\newline\url{https://www.fq.math.ca/Scanned/10-1/hoggatt2.pdf}

\bibitem[Vorobi02]{Vorobi02}Nicolai N. Vorobiev, \textit{Fibonacci numbers},
translated from the Russian by Mircea Martin, Springer 2002.

\bibitem[WooZei09]{1}Philip Matchett Wood, Doron Zeilberger, \textit{A
translation method for finding combinatorial bijections}, Annals of
Combinatorics \textbf{13} (2009), pp. 383--402. \newline\url{http://www.math.rutgers.edu/~zeilberg/mamarim/mamarimhtml/trans-method.html}
\end{thebibliography}
\end{document}